\documentclass[12pt]{amsart}
\usepackage{amsmath, amsthm, amscd, amsfonts}
\newtheorem{Lemma}{Lemma}[section]
\newtheorem{Theorem}{Theorem}[section]
\newtheorem{Corollary}{Corollary}[section]

 \newfont{\Bbbb}{msbm10 scaled\magstephalf}
 \def\be{\begin{equation}}
 \def\ee{\end{equation}}

\def\bn{\begin{equation}}
\def\en{\end{equation}}
\def\br{\begin{center}}
\def\er{\end{center}}
\def\by{\begin{array}}
\def\ey{\end{array}}

\def\begy{\begin{eqnarray}}
\def\endy{\end{eqnarray}}
\def\bey*{\begin{eqnarray*}}
\def\eny*{\end{eqnarray*}}
\def\ber{\begin{tabular}}
\def\enr{\end{tabular}}

\def\bt{\begin{flushright}}
\def\et{\end{flushright}}

 \def\bea{\begin{equation}\begin{array}{ll}}
 \def\eea{\end{array}\end{equation}}

\begin{document}
\title[Essential Norms of Weighted composition operators]{Essential Norms of Weighted Composition Operators between Hardy Spaces in the unit Ball}
\author[Z.S. Fang and Z.H.Zhou]{Zhong-Shan Fang \and Ze-Hua Zhou $^*$ }
\address{\newline Department of Mathematics\newline
Tianjin Polytechnic University
\newline Tianjin 300160\newline P.R. China.}

\email{fangzhongshan@yahoo.com.cn}

\address{\newline Department of Mathematics\newline
Tianjin University
\newline Tianjin 300072\newline P.R. China.}
\email{zehuazhou2003@yahoo.com.cn}

\keywords{Essential norm, Weighted composition operator, Bloch-type
space, Hardy space, Several complex variables}

\subjclass[2000]{Primary: 47B38; Secondary: 47B33,26A16, 32A16,
32A26, 32A30, 32A37, 32A38, 32H02.}

\date{}
\thanks{\noindent $^*$Ze-Hua Zhou, corresponding author. Supported in part by the National Natural Science Foundation of
China (Grand Nos.10671141, 10371091), and LiuHui Center for Applied
Mathematics, Nankai University \& Tianjin University.}

\begin{abstract}
 Let $\varphi(z)=(\varphi_1(z),\cdots,\varphi_n(z))$ be a
holomorphic self-map of $B_n$ and $\psi(z)$ a holomorphic function
on $B_n$, and $H(B_n)$ the class of all holomorphic functions on
$B_n$, where $B_n$ is the unit ball of $C^n$, the weight composition
operator $W_{\psi,\varphi}$ is defined by $W_{\psi,\varphi}=\psi
f(\varphi)$ for $f\in H(B_n)$. In this paper we estimate the
essential norm for the weighted composition operator
$W_{\psi,\varphi}$ acting from the Hardy space $H^p$ to $H^q$
($0<p,q\leq \infty$). When $p=\infty$ and $q=2$, we give an exact
formula for the essential norm. As their applications, we also
obtain some sufficient and necessary conditions for the bounded
weighted composition operator to be compact from $H^p$ to $H^q$.
\end{abstract}

\maketitle


\section{Introduction}

Let $B_n$ be the unit ball of $C^n$ with boundary $\partial B_n$,
$\sigma$ the normalized rotation invariant measure on $\partial
B_n$. The class of all holomorphic functions on domain $B_n$ will be
denoted by $H(B_n)$. Let
$\varphi(z)=(\varphi_1(z),\cdots,\varphi_n(z))$ be a holomorphic
self-map of $B_n$ and $\psi(z)$ is in $H(B_n)$. Multiplication
operator, Composition operator and weighted composition operator are
defined as follows:
$$M_{\psi}(f)(z)=\psi(z)\cdot f(z);$$
$$
C_\varphi(f)(z)=f(\varphi(z));$$
$$W_{\psi,\varphi}(f)(z)=\psi(z)\cdot f(\varphi(z))
$$
for any $ f\in H(B_n)$ and $z\in B_n$.

If let $\psi\equiv 1,$ then $W_{\psi,\phi}=C_{\phi};$  if let
$\phi=Id$, then $W_{\psi,\phi}=M_{\psi}.$ So we can regard weighted
composition operator as a generalization of a multiplication
operator and a composition operator. It is easy to show that
$C_{\phi}$ and $W_{\psi,\phi}$ take $H(B_n)$ into itself. Shapiro's
monograph \cite{Sh1} gives an interesting account of these
developments. See also Cowen and MacCluer's book \cite{Cow} for a
comprehensive treatment of these and other related problems with
composition operators.

In the recent years, boundedness and compactness of composition
operators between several spaces of holomorphic functions have been
studied by many authors: by Smith \cite{sm1}  between Bergman and
Hardy spaces, by Jarchow and Ried \cite{JR} between generalized
Bloch-type spaces and Hardy spaces, between Bloch spaces and Besov
spaces and BMOA and VMOA in Tian's thesis \cite{JR}, on BMOA by
Simth \cite{sm2}, and by Simth and Zhao \cite{sm3} from Bergman and
Hardy spaces and Bloch space into $Q_p$ spaces. All of papers above
focus on studying the composition operators in function spaces for
$1$-dimensional case.

More recently,  there have been many papers focused on studying the
same problems for $n$-dimensional case : by Luo and Shi \cite{LS1}
between Hardy spaces on the unit ball,\cite{LS2} weighted Bergman
spaces on bounded symmetric domains, by Zhou and
Shi\cite{zs1}\cite{zs2}\cite{zs3} on the Bloch space in polydisk or
classical symmetric domains, Gorkin and MacCluer \cite{GorM} between
hardy spaces in the unit ball, and Lipschitz space in polydisc by
Zhou \cite{z4}. In all these works the main goal is to relate
function theoretic properties of $\phi$ to boundedness and
compactness of $C_{\phi}$.

The essential norm of an operator T is by definition its distance to
the compact operators; that is $$ ||T||_e:=\inf\{||T-K||: K\;\;
compact\}.
$$
Notice that $\|T\|_e=0$ if and only if $T$ is compact, so that
estimates on $\|T\|_e$ lead to the conditions for $T$ to be compact.

In general, there is no easy way to determine the essential norms of
composition operator or weighted composition operator.

Let $f$ be in $H(B_n)$. For $0<p<\infty$, $f$ is said to be in the
Hardy space $H^p(B_n)$ provided that
$$
||f||^{p}_{p}=\sup\limits_{0<r<1}\int_{\partial B_n}|f(r\xi)|^p
d\sigma(\xi)<\infty .
$$
The Banach space of bounded holomorphic functions on $B_n$ in the
sup norm is donated by $H^{\infty}$.

When $f\in H^p$, then $f$ has radial limits at almost every
([d$\sigma$]) point of $\partial B_n$, and its $H^p$ norm is also
given by the $L^p(d\sigma)$ norm of its radial limit function
$f^{\ast}$. That is $$||f||^{p}_{p}=\int_{\partial
B_n}|f^{\ast}(\xi)|^p d\sigma(\xi).$$  Typically we continue to
write $f(\xi)$ for the radial limit; occasionally for clarity we use
the special notation $f^{\ast}(\xi)$ for $\lim_{r\rightarrow
1}f(r\xi)$. In the whole of paper, $E=\{\xi \in \partial
B_n:|\varphi^{\ast}(\xi)|=1 \}$, which we call it the extreme set of
$\varphi$.

It is well known that $C_{\varphi}$ is always bounded on $H^{p}(D)$
for $0<p\leq \infty$, this is a consequence of a theorem of J.
Littlewood, see \cite{Cow}, where $D=B_1$ is an unit disk. In 1987,
J.Shapiro \cite{Sh2} determined precisely when $C_{\varphi}$ acts
compactly on $H^p(D)$, for $p<\infty$, and gave a formula for the
essential norm of $C_{\varphi}$ acting on $H^2(D)$ in terms of the
Nevanlinna counting function for $\varphi$. In 2002, L. Zheng
\cite{Zhe} proved the essential norm of $C_{\varphi}$ acting on
$H^{\infty}(D)$ is $1$ whenever $C_{\varphi}$ is not compact on
$H^{\infty}(D)$ (equivalently, whenever $\|\varphi\|_{\infty}=1$);
it is also true when $D$ is replaced by the unit ball \cite{GMS}.
For $\infty\geq p>q>0,$ $C_{\varphi}$ acting from $H^p(D)$ to
$H^q(D)$ will of course be bounded. H. Jarchow \cite{JR} and T.
Goebeler \cite{Goe} shew independently that $C_{\varphi}$ is compact
if and only if $|E|=0$.

It seems reasonable to expect the essential norm to be given by a
formula that involves $|E|$. In fact, P.Gorkin and B.MacCluer
\cite{GorM} pointed out the essential norm of $C_{\varphi}$ acting
from $H^{\infty}(D)$ to $H^2(D)$ is precisely $|E|^{\frac{1}{2}}$,
and they have obtained the same results in the setting of Hardy
spaces $H^p(B_n)$ (we write it $H^p$ in the following) and also gave
some simple estimates for the essential norm of a composition
operator acting from $H^{\infty}$ to $H^q$ for $q\neq 2$ and for
$q<p<\infty,$ from $H^p$ to $H^q$ under a natural additional
condition. Here the additional condition is that there exists
$0<p<\infty$ such that $C_\varphi: H^p\rightarrow H^p$ is bounded,
which is naturally satisfied in the case $n=1$. This assumption has
two properties of interest to us:

(1) No set of positive measure in $\partial B_n$ is mapped by
$\varphi^{\ast}$ to a set of measure 0 in $\partial B_n$ (see
Corollary 3.38 of \cite{GorM});

(2) If $f\in H^p(B_n)$,then for a.e. $[d\sigma] \xi \in \partial
B_n$,
 $(f\circ\varphi)^{\ast}(\xi)=f^{\ast}(\varphi^{\ast}(\xi))$ (see
Lemma 1.6 in \cite{Mac}).

In our paper, in addition to extend corresponding cases in
\cite{GorM} to the weighted composition operator, we also get the
lower estimates for the essential norm of a weighted composition
operator from $H^p$ to $H^q$ for $1<p\leq q\leq \infty$.

The remainder of the present paper is assembled as follows: In
section 2, we refer the reader some Lemmas which needs in next
sections. In section 3, we will show that the essential norm of the
bounded weighted composition operator $W_{\psi,\varphi}$ is
precisely $(\mu_{\psi,\varphi,2}(\varphi(E)))^{1/2}$ for the case
$p=\infty, q=2$ (Theorem 3.1), and give a estimate for the case
$p=\infty, q\neq 2$ (Theorem 3.2). In section 4, we give the upper
estimate for the case $1<p<\infty$ (Theorem 4.1) and lower estimate
for the case $1<q<p<\infty$ (Theorem 4.2). The fundamental ideas of
the proof are those used by Gorkin and MacCluerin in \cite{GorM},
but some new techniques are still used in this section because of
the citation of the new measure induced by $\psi$ and $\phi$ and the
difference between weighted composition operator and composition
operator. If $\psi=1$, $W_{\psi,\varphi}=C_{\varphi}$, we can
completely the corresponding results in \cite{GorM}.

In sections 5 and 6 (not be considered in Gorkin and MacCluerin's
paper), using different methods, we obtain some estimates for the
essential norms of the weighted operator acting from $H^p$ to
$H^{\infty}$ for $p>1$ (Theorem 5.2) and from $H^p$ to $H^q$ for
$1<p\leq q<\infty$ (Theorem 6.2).

All of them are done under the same additional condition. As their
applications, we also obtained some sufficient and necessary
conditions for the weighted composition operator to be compact from
$H^p$ to $H^q$ for the above cases. For convenience, we always
abbreviate $H^p(B_n)$ to $H^p$.

\section{Some Lemmas}

\begin{Lemma} Let $\varphi$ is holomorphic self-map of
$B_n$ and  $\psi \in H^p$, where $0<p<\infty$. For any measurable
subset $E$ of $\partial B_n$, denote $\mu_{\psi,\varphi,p}(E) =
\int_{\varphi^{-1}(E)\cap
\partial B_n}|\psi|^{p}d\sigma$. Then
 $$\int_{\overline{B}_n}gd\mu_{\psi,\varphi,p}=\int_{\partial
B_n}|\psi|^{p}(g\circ\varphi)d\sigma,$$ where $g$ is an arbitrary
measurable positive function in $\overline{B}_n$.\end{Lemma}

{\bf Proof}\hspace*{4mm} If $g$ is a measurable simple function
defined on $\overline{B}_n$ given by
$g=\sum\limits_{i=1}^{n}\alpha_{i}\chi_{E_{i}}$, then
\begin{eqnarray*}
\int_{\overline{B}_n}gd\mu_{\psi,\varphi,p}
&=&\sum_{i=1}^{n}\alpha_{i}\mu_{\psi,\varphi,p}(E_{i}) =
\sum_{i=1}^{n}\alpha_{i}\int_{\varphi^{-1}(E_{i})\cap \partial
B_n}|\psi|^{p}d\sigma\\
&=&\int_{\partial
B_n}|\psi|^{p}(\sum_{i=1}^{n}\alpha_{i}\chi_{\varphi^{-1}(E_{i})\cap\partial
B_n})d\sigma\\
&=& \int_{\partial B_n}|\psi|^{p}(g\circ\varphi)d\sigma .
\end{eqnarray*}
Now, if $g$ is a measurable positive function in $\overline{B}_n$,
then we can take an increasing sequence $\{g_{m}\}$ of positive and
simple functions such that $g_{m}(z) \rightarrow g(z) $ for all $z
\in \overline{B}_n$, it follows that
$$\int_{\overline{B}_n}g_{m}d\mu_{\psi,\varphi,p} \rightarrow
\int_{\overline{B}_n}gd\mu_{\psi,\varphi,p}.$$ On the other hand,
$|\psi|^{p}(g_{m}\circ\varphi)$ is an increasing sequence such that
$$|\psi(z)|^{p}(g_{m}(\varphi(z)) \rightarrow
|\psi(z)|^{p}(g(\varphi(z)) $$ for all $z \in \overline{B}_n$, so
$$\int_{\overline{B}_n}g_{m}d\mu_{\psi,\varphi,p}=\int_{\partial
B_n}|\psi|^{p}(g_{m}\circ\varphi)d\sigma \rightarrow \int_{\partial
B_n}|\psi|^{p}(g\circ\varphi)d\sigma.$$ And the conclusion follows
by the uniqueness of the limit.

\begin{Lemma} (See p116 in \cite{Zhu2}) Suppose $0<p<\infty$ and $f\in H^p.$
Then $|f(z)|\leq \frac{||f||_p}{(1-|z|^2)^{n/p}}$ for all $z\in
B_n$.\end{Lemma}

\begin{Lemma} Let $\Omega$ be a domain in $C^n$, $f\in H(\Omega)$. If a
compact set $K$ and its neighborhood $G$ satisfy $K\subset
 G\subset\subset\Omega$ and $\rho=dist(K,\partial G)>0$, then
$$\sup\limits_{z\in K}|\frac{\partial
f}{\partial z_{j}}(z)|\leq\frac{\sqrt{n}}{\rho}\sup\limits_{z\in
G}|f(z)|.$$\end{Lemma}

{\bf Proof}\hspace{4mm} Since $\rho=dist(K,\partial G)>0,$ for any
$a\in K$, the polydisc
$$
P_a=\left\{(z_1, \cdots, z_n)\in C^n: |z_j-a_j|
<\displaystyle\frac{\rho}{\sqrt{n}}, j=1,\cdots,n\right\}
$$
is contained in $G$. Using Cauchy inequality, we have
$$
\left|\displaystyle\frac{\partial f}{\partial z_j}(a)
\right|\leq\displaystyle\frac{\sqrt{n}}{\rho}
\sup\limits_{z\in\partial_0 P_a}|f(z)|\leq
\displaystyle\frac{\sqrt{n}}{\rho}\sup\limits_{z\in G}|f(z)|.
$$
So the Lemma follows.

\begin{Lemma} For fixed $0<\delta<1$, let $G=\{z\in B_n: |z|\leq1-\delta\}$. Then
$$\lim\limits_{r\rightarrow 1}\sup\limits_{z\in G}|f(z)-f(rz)|=0$$
for any $f\in H^p(B_n)$.\end{Lemma}

{\bf Proof}
\begin{eqnarray*}
\sup\limits_{z\in G}|f(z)-f(rz)|&=&\sup\limits_{z\in
G}|\sum\limits^{n}_{j=1}(f(rz_1,rz_2,\cdots,rz_{j-1},z_{j},\cdots,z_n)\\
&-&f(rz_1,rz_2,\cdots,rz_{j},z_{j+1},\cdots,z_n))|\\
&\leq& \sup\limits_{z\in
G}\sum\limits^{n}_{j=1}|\int^{1}_{r}|z_{j}\frac{\partial f}{\partial
z_{j}}(rz_1,rz_{j-1},tz_j,z_{j+1},\cdots,z^n)dt|\\
&\leq& (1-r)n\sup\limits_{z\in G}|\frac{\partial f}{\partial
z_{j}}(z)|.
\end{eqnarray*}
Define $G_1=\{z\in B_n:|z|\leq 1-\frac{\delta}{2}\}$, then $G\subset
G_1$  and  $dist(G,\partial G_1)=\frac{\delta}{2}$.

It follows from Lemma 2.3 that
$$\sup\limits_{z\in G}|\frac{\partial f}{\partial z_{j}}(z)| \leq
\frac{2\sqrt{n}}{\delta}\sup\limits_{z\in G_1}|f(z)|.$$ If
$p=\infty$, then
$$\sup\limits_{z\in
G}|f(z)-f(rz)|\leq\frac{2(1-r)n\sqrt{n}}{\delta}||f||_{\infty}.$$
For $0<p<\infty$, it follows from Lemma 2.2 that
\begin{eqnarray*}
\sup\limits_{z\in
G}|f(z)-f(rz)|&\leq&\frac{2(1-r)n\sqrt{n}}{\delta}\sup\limits_{z\in
G_1}\frac{||f||_{p}}{(1-|z|^2)^{n/p}}\\
&\leq&\frac{2(1-r)n\sqrt{n}}{\delta}\sup\limits_{z\in
G_1}\frac{||f||_{p}}{(1-|z|)^{n/p}}\\
&\leq&\frac{2(1-r)n\sqrt{n}}{\delta}\frac{||f||_{p}}{(\frac{\delta}{2})^{n/p}}.
\end{eqnarray*}
Let $r\rightarrow 1$, the conclusion follows.

\begin{Lemma} (See corollary 1.3 in \cite{Cow}) A sequence in a
reflexive functional Banach space converges weakly if and only if it
is bounded and converges point-wise.\end{Lemma}

\begin{Lemma} Assume $\{f_m\}$ is a bounded sequence in $H^p(B_n)
(p>1)$, and $\{f_m\}$ converges weakly to $0$, then for any compact
operator $K$ from $H^p(B_n)$ to $Y$ ($Y$ is a normalized linear
space), we have $\|Kf_m\|_Y\rightarrow 0$.\end{Lemma}

{\bf Proof}\hspace*{4mm}This is easily followed by Lemma 2.5 and the
property of compact operator.

\section{From $H^{\infty}$ to $H^q$}

\textbf{Case 1. $p=\infty, q=2$}

It is well known that for any  $f\in H(B_n)$, $f$ has homogeneous
expansion $f(z)=\sum\limits^\infty_{s=0}F_s(z),$ where $F_s(z)$ is
the homogeneous polynomial
$\sum\limits_{|\alpha|=s}c(\alpha)z^{\alpha},$
$z^{\alpha}=z^{\alpha_1}\cdots z^{\alpha_n},$
$\alpha=(\alpha_1,\cdots,\alpha_n),$ and
$|\alpha|=\alpha_1+\cdots\alpha_n.$

If $f\in H^2(B_n),$ then
$$||f||^2_2=\sum_{\alpha}|c(\alpha)|^2||z^{\alpha}||^2_2,$$ where
$$||z^{\alpha}||^2_2=\frac{(n-1)!\alpha!}{(n-1+|\alpha|)!},$$ where
$\{\frac{z^{\alpha}}{||z^{\alpha}||_2}\}$ is an orthonormal basis
for $H^2(B_n)$, and $c(\alpha)=D^{\alpha}f(0)/\alpha!$ with $
\alpha!=\alpha_1!\cdots \alpha_{n}!$. If necessary, we refer the
reader to see \cite{Rud}.

For $m$ a positive integer, define the operators from $H^2(B_n)$ to
itself:
$$R_m(\sum^{\infty}_{s=0}F_s)=\sum^{\infty}_{s=m+1}F_s$$ and $$Q_m=I-R_m.$$ It is easy to show that $R_m$ is compact and
$\|R_m\|=1$.

\begin{Lemma} $W_{\psi,\varphi}: H^{\infty}\rightarrow H^{q},$
$0<q<\infty$ is bounded if and only if $\psi \in H^q$.\end{Lemma}

{\bf Proof}\hspace*{4mm}If $W_{\psi,\varphi}$ is bounded, let $f=1$,
then $W_{\psi,\varphi}f=\psi f(\varphi)=\psi\in H^q$. Conversely,
apparently we have $||W_{\psi,\varphi}f||_q\leq
||\psi||_q||f||_{\infty}$ for any $f\in H^{\infty}$, that is,
$||W_{\psi,\varphi}||\leq||\psi||_q$.

Using the same methods as that of Gorkin-MacCluer in \cite{GorM},
with minor modifications, we can obtain the following Lemmas 3.2 and
3.3. But for the reader's convenience, we give still the detail
proof for the results.
\begin{Lemma} If $W_{\psi,\varphi}: H^{\infty}\rightarrow H^{2}$ and
$\psi \in H^2$, then
$$\|W_{\psi,\varphi}\|_{e}=\lim\limits_{m\rightarrow \infty
}\|R_mW_{\psi,\varphi}\|.$$\end{Lemma}

{\bf Proof}\hspace*{4mm} On one hand, by hypothesis and Lemma 3.1,
we know $W_{\psi,\varphi}$ is bounded, so the compactness of $Q_m$
implies that $Q_mW_{\psi,\varphi}$ is also compact,
$$\|W_{\psi,\varphi}\|_{e}=\|(R_m+Q_m)W_{\psi,\varphi}\|_{e}=\|R_mW_{\psi,\varphi}\|_{e}
\leq\|R_mW_{\psi,\varphi}\|,$$ it follows that
$$\|W_{\psi,\varphi}\|_{e}\leq \liminf_{n\rightarrow \infty
}\|R_mW_{\psi,\varphi}\|.$$

On the other hand, let $K: H^{\infty}\rightarrow H^2$ be compact.
Since $\|R_m\|=1$,
\begin{eqnarray*}&&||W_{\psi,\varphi}-K||\geq||R_m(W_{\psi,\varphi}-K)||\\
&&=||R_mW_{\psi,\varphi}-R_mK||
\geq||R_mW_{\psi,\varphi}||-||R_mK||.\end{eqnarray*} Note that $K$
is compact, the image of the unit ball in $H^{\infty}$ under $K$ has
compact closure in $H^2$. Since $||R_m||=1$ and $R_mK$ tends to $0$
point-wise in $H^2$, $R_mK$ tends to $0$ uniformly on the unit ball
of $H^{\infty}$, that is $||R_mK|| \rightarrow 0$ as $n\rightarrow
\infty$. It follows that $$||W_{\psi,\varphi}||_{e}\geq
\limsup_{m\rightarrow \infty }||R_mW_{\psi,\varphi}||,$$ this
completes the proof.

\begin{Lemma} For $W_{\psi,\varphi}: H^{\infty}\rightarrow H^{2}$ and
$\psi \in H^2$, if $k$ is fixed positive integer and $g$ is any
non-constant holomorphic function on $B_n$ with $||g||_{\infty}\leq
1$, then $||Q_kW_{\psi,\varphi}(g^m)||_{2}\rightarrow 0$ as
$m\rightarrow \infty$.\end{Lemma}

{\bf Proof}\hspace*{4mm} If $\alpha$ is a multi-index with
$|\alpha|\leq k$, then
$$||z^{\alpha}||_{2}^{2}=\frac{(n-1)! \alpha
!}{(n-1+|\alpha|)!}\leq(k!)^{n}\equiv c(n,k).$$
 Since $\overline{D}^n(0,\frac{1}{2n}) \subseteq B_n$ and Cauchy's
estimates, for any holomorphic function $F$ in $B_n$, we have
$$\frac{D^{\alpha}F(0)}{\alpha !}\leq
(2n)^{|\alpha|}||F||_{\infty,\overline{D}^n(0,\frac{1}{2n})}$$
where$||F||_{\infty,\overline{D}^n(0,\frac{1}{2n})}$ denotes the
maximum modulus of $F$ on the polydisc
$\overline{D}^n(0,\frac{1}{2n})$. Since the series coefficients for
$F$ are $c(\alpha)=\frac{D^{\alpha}F(0)}{\alpha !}$, we get the
series coefficients for $\psi\cdot g^m\circ\varphi$ are bounded
above by $$(2n)^{|\alpha|}||\psi\cdot
g^m\circ\varphi||_{\infty,\overline{D}^n(0,\frac{1}{2n})}.$$

Let $c=max|\psi|$ and $s=max|g\circ\varphi|$ on
$\overline{D}^n(0,\frac{1}{2n})$, then $s<1$ by hypothesis. This
implies that $||\psi\cdot
g^m\circ\varphi||_{\infty,\overline{D}^n(0,\frac{1}{2n})}\leq cs^m$,
which tends to $0$ as $m\rightarrow \infty$. For fixed $k$,
$||Q_kW_{\psi,\varphi}(g^m)||_{2}^{2}=\sum_{|\alpha|\leq
k}|c(\alpha)|||z^{\alpha}||_{2}^{2}$, where $c(\alpha)$ is the
coefficients of $z^{\alpha}$ in the expansion of
$\psi\cdot(g\circ\varphi)^m$. By the above estimate, we have
$$||Q_kW_{\psi,\varphi}(g^m)||_{2}^{2}\leq \sum_{|\alpha|\leq
k}((2n)^kcs^m)^2c(n,k)\leq c'(n,k)s^{2m}.$$
 For fixed $k$, the last
expression tends to $0$ as $m\rightarrow \infty$.

\begin{Lemma}Let $\epsilon>0$, set $E_{\epsilon}=\{\xi \in
\partial
B_n:|\varphi(\xi)|\geq1-\epsilon \}$ and let $E_{\epsilon}^{c}$
denote its complement in $\partial B_n$, $\psi \in H^2 $. Define an
operator $K: H^{\infty}\rightarrow H^{2}$ by
$K(f)=P(\chi_{E_{\epsilon}^{c}}\psi\cdot(f\circ\varphi))$, where $P$
is the orthogonal projection of $L^2$ onto $H^2$ (where we identify
a function in $H^2$ with its radial limit function). Then $K$ is
compact from $H^{p}$ to $H^2$, for any $2<p\leq\infty$.\end{Lemma}

{\bf Proof}\hspace*{4mm} Let $\{f_m\}$ be a sequence from the unit
ball of $H^p$. By Lemma 2.4, $\{f_m\}$ is a normal family when
$2<p<\infty$, and this is obviously true for $p=\infty$. So there is
a subsequence which converges uniformly on compact subset of $B_n$,
to say $f$. For simplicity we still denote this subsequence as
$\{f_m\}$. Clearly $f\in H^p$. So
\begin{eqnarray*}
 ||Kf_m-Kf||_{2}^{2} &\leq &
||P||^2||\chi_{E_{\epsilon}^{c}}\psi\cdot(f_m\circ\varphi)-\chi_{E_{\epsilon}^{c}}\psi\cdot(f\circ\varphi)||_{2}^{2}\\
 & \leq & \int_{\partial
B_n}|\chi_{E_{\epsilon}^{c}}\psi\cdot(f_m\circ\varphi)-\chi_{E_{\epsilon}^{c}}\psi\cdot(f\circ\varphi)|^2d
\sigma\\
&=&\int_{E_{\epsilon}^{c}}|\psi\cdot(f_m\circ\varphi)-\psi\cdot(f\circ\varphi)|^2d\sigma.
\end{eqnarray*}
Since $\{f_m\}$ are uniformly bounded on $E_{\epsilon}^{c}$ and
$\psi \in H^2$, the above expression tends to $0$ as $n\rightarrow
\infty$ by Lebesgue's dominated convergence theorem. This verifies
the compactness of $K$.

\begin{Theorem} For $W_{\psi,\varphi}: H^{\infty}\rightarrow H^{2}$
and $\psi \in H^2$, then
$||W_{\psi,\varphi}||_{e}=(\mu_{\psi,\varphi,2}(\varphi(E)))^{1/2}$,
where $E=\{\xi \in \partial B_n:|\varphi^{\ast}(\xi)|=1
\}.$\end{Theorem}

{\bf Proof}\hspace*{4mm}we consider the lower estimate first.

Let $g$ be a non-constant inner function on $B_n$ and set $h=g^m$
for a positive integer $m,$ then
\begin{eqnarray*}
||W_{\psi,\varphi}(g^m)||_{2}^{2} &=&  \int_{\partial
B_n}|\psi^{\ast}\cdot(h^{\ast}\circ\varphi^{\ast})|^2d\sigma=
\int_{\overline{B}_n}|h^{\ast}|d\mu_{\psi,\varphi,2}\\
&\geq &\int_{\varphi(E)}|h^{\ast}|d\mu_{\psi,\varphi,2}\geq
\mu_{\psi,\varphi,2}(\varphi(E))
\end{eqnarray*}
where the last inequality follows by the fact that $|h^{\ast}|=1$
a.e $[d\mu]$ on $\varphi(E)$, this is true that $h$ is inner and the
restriction of $\mu_{\psi,\varphi,2}$ to $\partial B_n$ is
absolutely continuous
 with respect to $\sigma$.

In fact, for any measurable subset $E$ of $\partial B_n$,
$$\mu_{\psi,\varphi,2}(E)=\int_{\varphi^{-1}(E)\cap
\partial
 B_n}|\psi|^2d\sigma,$$ by hypothesis of
 $C_{\varphi}$, if $\sigma(E)=0$, then $\sigma(\varphi^{-1}(E))=0$, and
$\mu_{\psi,\varphi,2}(E)=0$ follows. So
\begin{eqnarray*}||R_kW_{\psi,\varphi}||&\geq&||R_kW_{\psi,\varphi}(g^m)||\geq||W_{\psi,\varphi}(g^m)||-||Q_kW_{\psi,\varphi}(g^m)||\\
 &\geq&
 \mu_{\psi,\varphi,2}(\varphi(E))-||Q_kW_{\psi,\varphi}(g^m)||.\end{eqnarray*}
for all $m$.

Fix $k$ and let $m\rightarrow \infty$ and apply Lemma 3.2 we obtain
$$||R_kW_{\psi,\varphi}||\geq(\mu_{\psi,\varphi,2}(\varphi(E)))^{1/2}$$
for any $k.$ Now let $k\rightarrow \infty,$ by Lemma 3.1 we have the
desired lower estimate on $||W_{\psi,\varphi}||_{e}.$

Now we turn to the upper estimate.

Take $K$ as in Lemma 3.3, for any $g\in H^{\infty}$ with
$||g||_{\infty}=1$, we have
\begin{eqnarray*}
||W_{\psi,\varphi}(g)-K(g)||_{2}&=& ||\psi\cdot
g\circ\varphi-P(\chi_{E_{\epsilon}^{c}}\psi\cdot(f\circ\varphi))||_2\\
& =&||P(\chi_{E_{\epsilon}}\psi\cdot(f\circ\varphi))||_2 \leq
||\chi_{E_{\epsilon}^{c}}\psi\cdot(f\circ\varphi)||_2\\
&=&(\int_{E_{\epsilon}}|\psi\cdot
g\circ\varphi|^2d\sigma)^{\frac{1}{2}}
\leq||g\circ\varphi||_{\infty}(\int_{E_{\epsilon}}|\psi|^2d\sigma)^{\frac{1}{2}}\\
&\leq&||g\circ\varphi||_{\infty}(\int_{\varphi^{-1}(\varphi(E_{\epsilon}))\cap
\partial B_n}|\psi|^2d\sigma)^{\frac{1}{2}}\\
&=&||g\circ\varphi||_{\infty}(\mu_{\psi,\varphi,2}(\varphi(E_{\epsilon})))^{1/2}.
\end{eqnarray*}
Let $\epsilon_{m}\downarrow 0$ and $K_{m}$ the corresponding
operator defined by
$$K_{m}(f)=P(\chi_{E_{\epsilon_{m}}^{c}}\psi\cdot(f\circ\varphi)).$$
For $p=\infty$ we have
$$||W_{\psi,\varphi}||_{e}\leq||W_{\psi,\varphi}-K_{m}||\leq
(\mu_{\psi,\varphi,2}(\varphi(E_{\epsilon_{m}})))^{1/2}$$ for all
$m,$ and let $m\rightarrow \infty$, as desired.

\begin{Corollary} $W_{\psi,\varphi}: H^{\infty}\rightarrow H^2$ is
compact if and only if $\psi\in H^2$ and
$\sigma(E)=0$.\end{Corollary}

{\bf Proof}\hspace*{4mm}If $W_{\psi,\varphi}$ is compact, it is
obviously bounded, it follows from Lemma 3.1 that $\psi\in H^2$.
From Theorem 3.1, the compactness of $W_{\psi,\varphi}$ implies
$\mu_{\psi,\varphi,2}(\varphi(E))=0$, so
$\sigma(\varphi^{-1}(\varphi(E))\cap
\partial B_n)=0$ (see 5.5.9 in \cite{Rud}), therefore $0\leq\sigma(E)\leq
\sigma(\varphi^{-1}(\varphi(E))\cap \partial B_n)=0$, $\sigma(E)=0$.

On the other hand, if $\psi\in H^2$, from the proof of theorem 3.1,
it follows that
$$||W_{\psi,\varphi}||_{e}\leq(\int_{E_{\epsilon}}|\psi|^2d\sigma)^{\frac{1}{2}}$$
when $\epsilon \rightarrow 0$ and since $\sigma(E)=0$, we get
$||W_{\psi,\varphi}||_{e}=0$, so $W_{\psi,\varphi}$ is compact.

In the above proof, set $\psi=1\in H^2$,  then
$||W_{1,\varphi}||_e=||C_\varphi||_e\leq \sigma(E)^{1/2}$. And if
set $\psi=1$ in theorem 3.1, then $$||C_\varphi||_e \geq
(\mu_{1,\varphi,2}(E))^{1/2}=\sigma(\varphi^{-1}(\varphi(E)))^{1/2},$$
so $\sigma(\varphi^{-1}(\varphi(E)))=\sigma(E)$, we have the
following Corollary

\begin{Corollary}(Theorem 1 in\cite{GorM}) $C_\varphi: H^{\infty}\rightarrow H^2$ is
bounded and
$$||C_\varphi||_e=\sigma(E)^{1/2}.$$\end{Corollary}

\textbf{Case 2. $p=\infty, q\neq 2$}

\begin{Theorem} Suppose $W_{\psi,\varphi}: H^{\infty}\rightarrow H^q
$ ($q>1$), and $\psi \in H^q$, then
$$
\frac{1}{2}(\mu_{\psi,\varphi,q}(\varphi(E)))^{1/q}\leq||W_{\psi,\varphi}||_{e}\leq2(\mu_{\psi,\varphi,q}
(\varphi(E)))^{1/q}.$$\end{Theorem}

{\bf Proof}\hspace*{4mm} We consider upper estimate first. Obviously
$W_{\psi,r\varphi}$ is compact for any fixed $0<r<1$. Let
$E_{\epsilon}=\{\xi \in
\partial B_n:|\varphi(\xi)|\geq1-\epsilon \}$ and let
$E_{\epsilon}^{c}$ denote its complement in $\partial B_n$. So
\begin{eqnarray*}
 ||W_{\psi,\varphi}-W_{\psi,r\varphi}||
 &=&
\sup\limits_{||f||_{\infty}=1}||(W_{\psi,\varphi}-W_{\psi,r\varphi})f||_q \\
&=&  \sup\limits_{||f||_{\infty}=1} (\int_{\partial B_n}|\psi(f\circ
\varphi)-\psi(f\circ (r\varphi))|^q d \sigma)^{1/q}
\\
&=& \sup\limits_{||f||_{\infty}=1}(\int_{E_{\epsilon}}|\psi(f\circ
\varphi)-\psi(f\circ (r\varphi))|^q d\sigma)^{1/q}\\
 &+&
\sup\limits_{||f||_{\infty}=1}(\int_{E_{\epsilon}^{c}}|\psi(f\circ
\varphi)-\psi(f\circ (r\varphi))|^q d\sigma)^{1/q}.
\end{eqnarray*}
Apply Lemma 2.4, we can choose $r$ sufficiently close to $1$ to make
the second term less than $\epsilon ||\varphi||_q$. For the first
term, the triangle inequality yields $$|f\circ
\varphi(\xi)-f\circ(r\varphi)(\xi)|\leq 2$$ So, the first term is
less than
$$2(\int_{E_{\epsilon}}|\psi|^q d \sigma)^{1/q}\leq
2(\int_{\varphi^{-1}(\varphi(E_{\epsilon}))\cap\partial B_n}|\psi|^q
d
\sigma)^{1/q}=2(\mu_{\psi,\varphi,q}(\varphi(E_{\epsilon})))^{1/q}.$$
Let $\epsilon_{m}\downarrow 0$, and $E_{\epsilon_m}=\{\xi \in
\partial B_n:|\varphi(\xi)|\geq1-\epsilon_m
\}$, then $\mu_{\psi,\varphi,q}(\varphi(E_{\epsilon_m})) \rightarrow
\mu_{\psi,\varphi,q}(\varphi(E))$, the upper estimate follows.

Now we turn to lower estimate. Let $f$ be a non-constant inner
function in $B_n$, $K$ is any compact operator. For any positive
integer $m$, the sequence $\{f^m\}$ are in
 the unit ball of $H^{\infty}$, So there exists a subsequence
$\{f^{m_k}\}$ such that $\{K(f^{m_k})\}$ converges in norm.
Therefore, given $\epsilon>0$, there exists $M$ such that
$||K(f^{m_k})-K(f^{m_l})||_q<\epsilon$ for any $k,l>M$. Fix $k>M$,
there exists $r$ with $0<r<1$ such that
 $(\psi(f\circ\varphi)^{m_k})_r(z)=\psi(rz)(f\circ\varphi(rz))^{m_k}$
satisfies
$$||(\psi(f\circ\varphi)^{m_k})_r||_q \geq
||(\psi(f\circ\varphi)^{m_k})||-\epsilon.$$ Thus, for $m\geq M$
\begin{eqnarray*}
||W_{\psi,\varphi}-K||&\geq&||(W_{\psi,\varphi}-K)\frac{
f^{m_k}-f^{m_l}}{2}||_q\\
& \geq &
(1/2)||(\psi(f\circ\varphi)^{m_k})-(\psi(f\circ\varphi)^{m_l})||_q-\epsilon/2\\
& \geq
&(1/2)(||(\psi(f\circ\varphi)^{m_k})||_q-||(\psi(f\circ\varphi)^{m_l})||_q)-\epsilon/2\\
& \geq &
(1/2)(||(\psi(f\circ\varphi)^{m_k})||_q-||(\psi(f\circ\varphi)^{m_l})_r||_q)-\epsilon.
\end{eqnarray*}
letting  $l\rightarrow \infty$  and $h=f^{m_k}$, we have
\begin{eqnarray*}
||W_{\psi,\varphi}-K|| &\geq&
(1/2)(||(\psi(f\circ\varphi)^{m_k})||_q-\epsilon\\
&=&  (1/2)(\int_{\partial
B_n}|\psi^{\ast}\cdot(h^{\ast}\circ\varphi^{\ast})|^q
d\sigma)^{1/q}-\epsilon\\
&=&(1/2) (\int_{\overline{B}_n}|h^{\ast}|^q
d\mu_{\psi,\varphi,q})^{1/q}-\epsilon\\
&\geq &(1/2)(\int_{\varphi(E)}|h^{\ast}|^q
d\mu_{\psi,\varphi,q})^{1/q}-\epsilon\\
&\geq&(1/2) (\mu_{\psi,\varphi,q}(\varphi(E)))^{1/q}-\epsilon
\end{eqnarray*}
  Now letting $\epsilon \rightarrow 0$ yields the result.

\begin{Corollary} $W_{\psi,\varphi}: H^{\infty}\rightarrow H^q
$ is compact if and only if $\psi\in H^q$ and
$\sigma(E)=0$.\end{Corollary}

{\bf Proof}\hspace*{2mm} Combining Lemma 3.1 and Theorem 3.2, the
corollary follows.

\begin{Corollary} (Theorems 2 and 3 \cite{GorM} $C_\varphi:
H^{\infty}\rightarrow H^q $ is bounded and
$$\frac{1}{2}\sigma(E)^{1/q}\leq||C_\varphi||_e\leq2\sigma(E)^{1/q}.$$\end{Corollary}

{\bf Proof}\hspace*{2mm} Let $\psi=1\in H^q,$ then
$W_{\psi,\varphi}=C_\varphi$, the corollary follows by Theorem 3.2.

\section{From $H^p$ to $H^q$ for $1<q<p<\infty$}

\begin{Theorem} Assume $W_{\psi,\varphi}: H^{p}\rightarrow H^q $
($1<p<\infty$) is bounded, then $||W_{\psi,\varphi}||_{e}\geq
(\mu_{\psi,\varphi,q}(\varphi(E)))^{1/q}$.\end{Theorem}

{\bf Proof}\hspace*{4mm}Let $g$ be a non-constant inner function on
$B_n$ and set $h=g^m$ for a positive integer $m$. Then $||g^m||_p=1$
for any $m,$ and $g^m$ converges weakly to $0$ as $m\rightarrow
\infty$, thus $||Kf_w||\rightarrow 0$ for any compact operator from
$H^P$ to $H^q$ when $|w|\rightarrow 1$. Like in Theorem 3.1, we have
\begin{eqnarray*}
||W_{\psi,\varphi}-K|| &\geq& \limsup\limits_{m\rightarrow
\infty}||(W_{\psi,\varphi}-K)(g^m)||_{q}\\
&\geq&\limsup\limits_{m\rightarrow
\infty}||W_{\psi,\varphi}(g^m)||_q-\limsup\limits_{m\rightarrow
\infty}||K(g^m)||_q \\
 &=& \limsup\limits_{m\rightarrow
\infty}||W_{\psi,\varphi}(g^m)||_q     =\limsup\limits_{m\rightarrow
\infty} (\int_{\partial
B_n}|\psi^{\ast}\cdot(h^{\ast}\circ\varphi^{\ast})|^qd\sigma)^{1/q}\\
&=& \limsup\limits_{m\rightarrow \infty}
(\int_{\overline{B}_n}|h^{\ast}|d\mu_{\psi,\varphi,q})^{1/q} \geq
\limsup\limits_{m\rightarrow
\infty} (\int_{\varphi(E)}|h^{\ast}|d\mu_{\psi,\varphi,q})^{1/q}\\
&\geq& (\mu_{\psi,\varphi,2}(\varphi(E)))^{1/q}.
\end{eqnarray*}This ends the proof.

\begin{Corollary} Assume $W_{\psi,\varphi}: H^p\rightarrow H^q, p>1,
0<q<\infty $ is compact, then $\sigma(E)=0$.\end{Corollary}

\textbf{Remark 1.} We will show that when  $0<p<q<\infty$ and
$W_{\psi,\varphi}: H^{p}\rightarrow H^q $ is bounded, then
$\mu_{\psi,\varphi,q}(\varphi(E))=0$ (see Corollary 6.1), So the
above estimate is useless.

\begin{Theorem} Suppose $1<q<p<\infty$ and there exists $r>q $ such
that $W_{\psi,\varphi}: H^{p}\rightarrow H^r $ ($1<p<\infty$) is
bounded, then $$||W_{\psi,\varphi}||_e\leq
||P||\cdot||W_{\psi,\varphi}||_{p,r}\cdot\sigma(E)^{\frac{r-q}{qr}}$$
where $P$ is the Szeg\"{o} projection of $L^q(\sigma)$ onto
$H^q$.\end{Theorem}

{\bf Proof}\hspace*{4mm}We consider the operator $K:
H^{p}\rightarrow H^{q}$ defined by
$$K(f)=P(\chi_{E_{\epsilon}^{c}}\psi\cdot(f\circ\varphi)),$$
where $P$ is the Szeg\"{o} projection of $L^q(\sigma)$ onto $H^q$.
Like in Lemma 3.3, $K$ is compact operator from $H^{p}$ to $H^{q}$.
So for any $g\in H^p$ with $||g||_p=1$, we have
\begin{eqnarray*}
||W_{\psi,\varphi}(g)-K(g)||_{q}&=& ||\psi\cdot
g\circ\varphi-P(\chi_{E_{\epsilon}^{c}}\psi\cdot(f\circ\varphi))||_q \\
&=&||P(\chi_{E_{\epsilon}}\psi\cdot(f\circ\varphi))||_q\\
&\leq&
||P||\cdot||\chi_{E_{\epsilon}^{c}}\psi\cdot(f\circ\varphi)||_q\\
&=&||P||\cdot(\int_{E_{\epsilon}}|\psi\cdot
g\circ\varphi|^qd\sigma)^{\frac{1}{q}}\\
&\leq&||P||\cdot(\int_{\partial
B_n}\chi_{E_{\epsilon}^{c}}|\psi\cdot
g\circ\varphi|^qd\sigma)^{\frac{1}{q}}\\
&\leq&||P||\cdot||W_{\psi,\varphi}(g)||_r\sigma(E_{\epsilon})^{\frac{r-q}{qr}}\\
&\leq&||P||\cdot||W_{\psi,\varphi}||_{p,r}\cdot\sigma(E_{\epsilon})^{\frac{r-q}{qr}}.
\end{eqnarray*}
Letting $\epsilon\rightarrow 0$ yields the conclusion.

\section{From $H^p$ to $H^{\infty}$}

\begin{Theorem} For $W_{\psi,\varphi}: H^p\rightarrow H^{\infty} $,
and $0<p<\infty$, then $W_{\psi,\varphi}$ is bounded if and only if
$\sup\limits_{z\in
B_n}\frac{|\psi(z)|}{(1-|\varphi(z)|^2)^{n/p}}<\infty$.\end{Theorem}

{\bf Proof}\hspace*{4mm} $"\Rightarrow"$ For any $w\in B_n$, define
$f_{w}(z)=\frac{(1-|w|^2)^{n/p}}{(1-<z,w>)^{2n/p}}$, and it is easy
to check $||f_{w}||_{p}=1$. So \begin{eqnarray*}C\geq
||W_{\psi,\varphi}||&=&\sup\limits_{||f||_{p}=1}||W_{\psi,\varphi}f||_{\infty}\geq\sup\limits_{z\in
B_n}||W_{\psi,\varphi}f_{w}||_{\infty}\\
&=&\sup\limits_{w\in B_n}\sup\limits_{z\in
B_n}|\psi(z)||f_{w}(\varphi(z))|\end{eqnarray*} setting
$w=\varphi(z)$, as desired.

$"\Leftarrow"$
\begin{eqnarray*}
||W_{\psi,\varphi}||&=&\sup\limits_{||f||_{p}=1}||W_{\psi,\varphi}||_{\infty}=\sup\limits_{||f||_{p}=1}\sup\limits_{z\in
B_n}|\psi(z)f(\varphi(z))|\\
&\leq& \sup\limits_{||f||_{p}=1}\sup\limits_{z\in
B_n}|\psi(z)|\frac{||f||_{p}}{(1-|\varphi(z)|^2)^{n/p}}
=\sup\limits_{z\in B_n}\frac{|\psi(z)|}{(1-|\varphi(z)|^2)^{n/p}}
\end{eqnarray*}

\begin{Theorem} For $W_{\psi,\varphi}: H^p\rightarrow H^{\infty} $
($p>1$), and $W_{\psi,\varphi}$ is bounded, then
\begin{eqnarray*}&&\lim\limits_{\delta \rightarrow
0}\sup\limits_{dist(\varphi(z),\partial
B_n)<\delta}\frac{|\psi(z)|}{(1-|\varphi(z)|^2)^{n/p}}
\leq||W_{\psi,\varphi}||_{e}\\
&&\hspace*{6mm}\leq 2\lim\limits_{\delta \rightarrow
0}\sup\limits_{dist(\varphi(z),\partial
B_n)<\delta}\frac{|\psi(z)|}{(1-|\varphi(z)|^2)^{n/p}}.\end{eqnarray*}\end{Theorem}

{\bf Proof}\hspace*{4mm}We consider the upper estimate first.

For any fixed $0<r<1$, it is easy to check that $W_{\psi,r\varphi}$
is compact. Thus $$||W_{\psi,\varphi}||_{e}\leq
||W_{\psi,\varphi}-W_{\psi,r\varphi}||.$$ Now for any $0<\delta<1$
\begin{eqnarray*}
||W_{\psi,\varphi}-W_{\psi,r\varphi}||&=&\sup\limits_{||f||_p=1}||(W_{\psi,\varphi}-W_{\psi,r\varphi})f||_{\infty}
\\
&=&\sup\limits_{||f||_p=1}\sup\limits_{z\in
B_n}|\psi(z)|\cdot|f(\varphi(z))-f(r\varphi(z))|\\
&\leq&
||\psi||_{\infty}\sup\limits_{||f||_p=1}\sup\limits_{dist(\varphi(z),\partial
B_n)\geq\delta}|f(\varphi(z))-f(r\varphi(z))|\\
&+& \sup\limits_{||f||_p=1}\sup\limits_{dist(\varphi(z),\partial
B_n)<\delta}|\psi(z)|\cdot|f(\varphi(z))-f(r\varphi(z))|.
     \end{eqnarray*}
From Lemma 2.4, we can choose $r$ sufficiently close to $1$ such
that the first term of the right hand side is less than any given
 $\epsilon$. And we denote the second term by $I$.
Then,
\begin{eqnarray*}
I&\leq&\sup\limits_{||f||_p=1}\sup\limits_{dist(\varphi(z),\partial
B_n)<\delta}|\psi(z)|\cdot(|f(\varphi(z))|+|f(r\varphi(z))|)\\
 &\leq&\sup\limits_{||f||_p=1}\sup\limits_{dist(\varphi(z),\partial
B_n)<\delta}|\psi(z)|(\frac{||f||_p}{(1-|\varphi(z)|^2)^{n/p}}+\frac{||f||_p}{(1-|r\varphi(z)|^2)^{n/p}})\\
 &\leq&2\sup\limits_{dist(\varphi(z),\partial
B_n)<\delta}\frac{|\psi(z)|}{(1-|\varphi(z)|^2)^{n/p}}.
\end{eqnarray*}
Now let $r\rightarrow 1$ first, then let $\delta \rightarrow 0$, we
get the desired upper estimate.

We now turn to the lower estimate.

Let $K$ be any compact operator from $H^p$ to $H^{\infty}$. For any
$w\in B_n$ define
$f_{w}(z)=\frac{(1-|w|^2)^{n/p}}{(1-<z,w>)^{2n/p}}$, it is easy to
check $||f_{w}||_{p}=1$ and $f_w$ converge weakly to $0$ as
$|w|\rightarrow 1$, thus $||Kf_w||\rightarrow 0$ when
$|w|\rightarrow 1$.

So for any $0<\delta<1$
\begin{eqnarray*}
||W_{\psi,\varphi}-K||&\geq& \limsup\limits_{|w|\rightarrow
1}||(W_{\psi,\varphi}-K)f_w||_{\infty}\\
&\geq&\limsup\limits_{|w|\rightarrow
1}||W_{\psi,\varphi}f_{w}||_{\infty}-\limsup\limits_{|w|\rightarrow
1}||Kf_{w}||_{\infty}\\
&=&\limsup\limits_{|w|\rightarrow 1}\sup\limits_{z\in
B_n}|\psi(z)||f_{w}(\varphi(z))|\\
&\geq&\limsup\limits_{|w|\rightarrow
1}\sup\limits_{dist(\varphi(z),\partial
B_n)<\delta}|\psi(z)||f_{w}(\varphi(z))|
\end{eqnarray*}
Let $\delta \rightarrow 0$ then $|\varphi(z)|\rightarrow 1$ and set
$w=\varphi(z)$, we obtain the lower estimate of
$||W_{\psi,\varphi}||_{e}$.

\begin{Corollary} Assume $W_{\psi,\varphi}: H^p\rightarrow
H^\infty $ is bounded, then it is compact if and only if $$
\lim\limits_{\delta \rightarrow
0}\sup\limits_{dist(\varphi(z),\partial
B_n)<\delta}\frac{|\psi(z)|}{(1-|\varphi(z)|^2)^{n/p}}=0.$$\end{Corollary}

\textbf{Remark 2.} If $||\varphi||_{\infty}<1$, then $E=\{z\in
\overline{B}_n| \varphi(z)=1\}=\emptyset$, without the loss of
generality, we set
$$\lim\limits_{\delta \rightarrow
0}\sup\limits_{dist(\varphi(z),\partial
B_n)<\delta}\frac{|\psi(z)|}{(1-|\varphi(z)|^2)^{n/p}}=0.$$

\section{From $H^p$ to $H^q$ for $1<p\leq q<\infty$}

\textbf{Definition} Let $\beta\geq 1$. A finite and positive measure
$\mu$ on is called a $\beta-Carleson$ measure. If there is a
constant $M\leq\infty$ such that $\mu(S_h(\xi))\geq Mh^{n\beta}$ for
all $\xi \in \partial B_n$ and $0<h<2$, and it is called vanishing
$\beta-Carleson$ measure if $\lim\limits_{h\rightarrow
0}\sup\limits_{\xi \in \partial
B_n}\frac{\mu(S_h(\xi))}{h^{n\beta}}=0.$

\begin{Lemma} (see corollary 2 in \cite{LS1}). Let $\mu$ be a
finite and positive measure on $\overline{B}_n$, and $0<p\leq
q<\infty$, then the following statement are equivalent:

(i)\hspace{2mm} $\mu$ is a bounded $\frac{q}{p}-Carleson$ measure.

(ii)\hspace{2mm} There is a constant $C<\infty$ so that
$$\int_{\overline{B}_n}|f|^p d\mu \leq C||f||_{p}^{q}$$
for all $f$ in $\overline{B}_n$.\end{Lemma}

\begin{Lemma} (\cite{X}) Suppose that $0< p\leq q<\infty$ and $W_{\psi,\varphi}:
H^p\rightarrow H^q$ is bounded, then the following conditions are
equivalent:

(i)\hspace{2mm}$\mu_{\psi,\varphi,q}$ is vanishing
$\frac{q}{p}-Carleson$ measure

(ii)\hspace{2mm}$W_{\psi,\varphi}: H^p\rightarrow H^q$ is compact
operator.\end{Lemma}

\begin{Theorem} For Fixed $0\leq p\leq q<\infty$, then the following
statement are equivalent:

(i)\hspace{2mm}$\mu_{\psi,\varphi,q}$ is a bounded
$\frac{q}{p}-Carleson$ measure.

(ii)\hspace{2mm}$W_{\psi,\varphi}: H^p\rightarrow H^q$ is bounded.

(iii) $$\sup\limits_{z \in
B_n}\int_{\overline{B}_n}\frac{(1-|z|^2)^{nq/p}}{|1-<w,z>|^{2nq/p}}d\mu_{\psi,\varphi,q}(w)<\infty.$$\end{Theorem}

{\bf Proof}\hspace*{4mm}$(i)\Rightarrow (ii)$ \\
From Lemma 3.4 , if $\mu_{\psi,\varphi,q}$ is bounded
$\frac{q}{p}-Carleson$ measure, then there exists constant $C$ such
that
$$\int_{\overline{B}_n}|f|^p d\mu_{\psi,\varphi,q} \leq C||f||_{p}^{q}$$
for any  $f \in H^p(B_n)$. Apply Lemma 2.1, and put $g=|f|^q$, we
have
$$\int_{\overline{B}_n}|f|^q d\mu_{\psi,\varphi,q} =\int_{\partial
B_n}|\psi|^{q}|f\circ\varphi|^q d\sigma=||W_{\psi,\varphi}f||_q^q.$$
So $$||W_{\psi,\varphi}(f)||_q\leq C^{1/q}||f||_p$$ for any $f \in
H^p(B_n)$. That is, $W_{\psi,\varphi}: H^p\rightarrow H^q$ is
bounded.

$(ii) \Rightarrow (iii)$ \\
For any $z\in B_n $ , set
$f_{z}(w)=\frac{(1-|z|^2)^{n/p}}{(1-<w,z>)^{2n/p}}$, then
$||f_{w}||_{p}=1$
\begin{eqnarray*}
C&\geq&
||W_{\psi,\varphi}||^q=\sup\limits_{||f||_{p}=1}||W_{\psi,\varphi}f||^{q}_{q}\geq\sup\limits_{z\in
B_n}||W_{\psi,\varphi}f_{z}||^{q}_{q}\\
&=&\sup\limits_{z\in B_n}(\int_{\partial
B_n}|\psi|^{p}|f_z\circ\varphi|^q d\sigma) =\sup\limits_{z\in
B_n}\int_{\overline{B}_n}|f_z|^q
d\mu_{\psi,\varphi,q}\\
&=&\sup\limits_{z \in
B_n}\int_{\overline{B}_n}\frac{(1-|z|^2)^{nq/p}}{|1-<w,z>|^{2nq/p}}d\mu_{\psi,\varphi,q}(w)
\end{eqnarray*}

$(iii) \Rightarrow (i)$ \\
Assume that $$M=\sup\limits_{z \in
B_n}\int_{\overline{B}_n}\frac{(1-|z|^2)^{nq/p}}{|1-<w,z>|^{2nq/p}}d\mu_{\psi,\varphi,q}(w)<\infty$$
we show that $\mu_{\psi,\varphi,q}$ is a bounded
$\frac{q}{p}-Carleson$ measure.

First let $z=0$, then $\mu_{\psi,\varphi,q}(\overline{B}_n)\leq M$.
Thus $\mu_{\psi,\varphi,q}$ is finite and hence
$\mu_{\psi,\varphi,q}(S_h(\xi))\leq M\leq 4Mh^{nq/p}$ for all $\xi
\in \partial B_n$ and $h\geq (\frac{1}{4})^{\frac{1}{nq/p}}$.
Suppose $h\leq (\frac{1}{4})^{\frac{1}{nq/p}}$ and $\xi \in \partial
B_n$. Let $\xi_0=(1-\frac{h}{2})\xi$, then for any $w\in S_h(\xi)$,
\begin{eqnarray*}
|1-<w,\xi_0>|&=&|1-\frac{h}{2}+\frac{h}{2}-<w,\xi_0>|\\
&=&|(1-\frac{h}{2})(1-<w,\xi>)+\frac{h}{2}|\\
&\leq& |(1-\frac{h}{2})h|+\frac{h}{2} \leq \frac{3h}{2}
\end{eqnarray*}
and  $ 1-|\xi_0|^2=(1-|\xi_0|)(1+|\xi_0|)\geq(1-|\xi_0|)$, we have
$$\frac{(1-|\xi_0|^2)^{nq/p}}{|1-<w,\xi_0>|^{2nq/p}}\geq
\frac{(1-|\xi_0|)^{nq/p}}{\frac{3h}{2}}=\frac{c}{h^{nq/p}}.$$ So
\begin{eqnarray*}
M&\geq&
\int_{\overline{B}_n}\frac{(1-|\xi_0|^2)^{nq/p}}{|1-<w,\xi_0>|^{2nq/p}}d\mu_{\psi,\varphi,q}(w)\\
&\geq&\int_{S_h(\xi)}\frac{c}{h^{nq/p}}d\mu_{\psi,\varphi,q} \geq
\frac{c\mu_{\psi,\varphi,q}(S_h(z))}{h^{nq/p}}.
\end{eqnarray*}
Therefore, $\mu_{\psi,\varphi,q}$ is bounded $\frac{q}{p}-Carleson$
measure.

\begin{Corollary} If $0<p<q<\infty$ and $W_{\psi,\varphi}:
H^{p}\rightarrow H^q $ is bounded, then
$\mu_{\psi,\varphi,q}(\varphi(E))=0.$\end{Corollary}

{\bf Proof}\hspace*{4mm}Denote $g$ the $Radon-Nikod\acute{y}m$
derivative of $\mu_{\psi,\varphi,q}|\partial B_n$ with respect to
$\sigma$, $\mu_{\psi,\varphi,q}$ is absolutely continuous with
respect to $\sigma$ on $\partial B_n$, so it follows that
$$g(b)=\lim\limits_{h\rightarrow
0}\frac{1}{\sigma(S_h(b))}\int_{S_h(b)}gd\sigma=\lim\limits_{h\rightarrow
0}\frac{\mu_{\psi,\varphi,q}(S_h(b))}{\sigma(S_h(b))}\geq\lim\limits_{h\rightarrow
0}Ch^{nq/p-n}=0$$ almost everywhere in $\partial B_n$. Where the
penultimate inequality uses the fact that $\sigma(S_h(b))$ is
roughly proportional to $h^n$(see P67 in \cite{Rud}). Now we have
$\mu_{\psi,\varphi,q}|\partial B_n$=0, the corollary is proved.

\begin{Theorem} For fixed $1< p\leq q<\infty$ and weighted
composition operator $W_{\psi,\varphi}: H^p\rightarrow H^q$ is
bounded, then$$||W_{\psi,\varphi}||_e\geq \lim\limits_{|w|
\rightarrow
1}\int_{\overline{B}_n}\frac{(1-|w|^2)^{nq/p}}{|1-<z,w>|^{2nq/p}}d\mu_{\psi,\varphi,q}(z).$$
\end{Theorem}

{\bf Proof}\hspace*{4mm}Let $K$ be any compact operator from $H^p$
to $H^{\infty}$. For any $w\in B_n$ define
$f_{w}(z)=\frac{(1-|w|^2)^{n/p}}{(1-<z,w>)^{2n/p}}$, it is easy to
check $||f_{w}||_{p}=1$ and $f_w$ converge weakly to $0$ as
$|w|\rightarrow 1$, thus $||Kf_w||\rightarrow 0$ when
$|w|\rightarrow 1$. So for any $0<\delta<1,$
\begin{eqnarray*}
||W_{\psi,\varphi}-K||&\geq& \limsup\limits_{|w|\rightarrow
1}||(W_{\psi,\varphi}-K)f_w||_q\\
&\geq&\limsup\limits_{|w|\rightarrow
1}||W_{\psi,\varphi}f_{w}||_q-\limsup\limits_{|w|\rightarrow
1}||Kf_{w}||_q\\
&=&\limsup\limits_{|w|\rightarrow 1}   \int_{\partial
B_n}|\psi(z)|^q\frac{(1-|w|^2)^{nq/p}}{|1-<\varphi(z),w>|^{2nq/p}}d\sigma(z)\\
&\geq&\limsup\limits_{|w|\rightarrow
1}\int_{\overline{B}_n}\frac{(1-|w|^2)^{nq/p}}{|1-<z,w>|^{2nq/p}}d\mu_{\psi,\varphi,q}(z)
\end{eqnarray*}
The conclusion follows.

We cannot give the upper estimate in the above form, but we have the
following theorem.

\begin{Theorem} Assume $1< p\leq q<\infty$ and $W_{\psi,\varphi}:
H^p\rightarrow H^q$ is bounded, then $W_{\psi,\varphi}:
H^p\rightarrow H^q$ is compact if and only if
$$\lim\limits_{|w|
\rightarrow1}\int_{\overline{B}_n}\frac{(1-|w|^2)^{nq/p}}{|1-<z,w>|^{2nq/p}}d\mu_{\psi,\varphi,q}(z)=0.$$\end{Theorem}

{\bf Proof}\hspace*{4mm}The necessary condition follows by theorem
6.2. We consider the sufficient condition. By Lemma 6.2, we only
have to show $\mu_{\psi,\varphi,q}$ is vanishing
$\frac{q}{p}-Carleson$ measure. From the proof of $(iii)\Rightarrow
(i)$ in theorem 6.1, for any $z\in \partial B_n$, set
$|z_0|=1-\frac{h}{2}$. Suppose
$$\lim\limits_{|w|
\rightarrow1}\int_{\overline{B}_n}\frac{(1-|w|^2)^{nq/p}}{|1-<z,w>|^{2nq/p}}d\mu_{\psi,\varphi,q}(z)=0.$$
That is , $\forall \epsilon>0,\exists 1>r>0,$ when $|w|>r$ we have
$$|\int_{\overline{B}_n}\frac{(1-|w|^2)^{nq/p}}{|1-<z,w>|^{2nq/p}}d\mu_{\psi,\varphi,q}(z)|<\epsilon.$$
When $h<2(1-r)$, for any $z\in \partial B_n$, the corresponding
$|z_0|>r$, so
\begin{eqnarray*}
\epsilon&>&\int_{\overline{B}_n}\frac{(1-|z_0|^2)^{nq/p}}{|1-<w,z_0>|^{2nq/p}}d\mu_{\psi,\varphi,q}(w)\\
&\geq&\int_{S_h(z)}\frac{c}{h^{nq/p}}d\mu_{\psi,\varphi,q}\\
&\geq& \frac{c\mu_{\psi,\varphi,q}(S_h(z))}{h^{nq/p}}.
\end{eqnarray*}
This is true for any $z\in \partial B_n$. So $\mu_{\psi,\varphi,q}$
is vanishing $\frac{q}{p}-Carleson$ measure.

\end{document}